\documentclass[10pt,twoside]{article}
\usepackage{graphicx}
\usepackage{amsmath}
\usepackage{Latex-document}
\def\volumeno{I}
\title{\bf Smoothed Analysis of Algorithms\vskip 6mm}
\author{{\bf Daniel A. Spielman}\thanks{Department of
Mathematics, MIT, USA. E-mail: spielman@math.mit.edu} \quad
Shang-Hua Teng\thanks{Department of Computer Science, Boston
University, USA.}\vspace*{-0.5cm}}
\date{\vspace{-8mm}}

\newdimen\pIR
\newcommand\StevesR{{\rm I\kern\pIR R}}

\newcommand{\qed}{\hspace*{\fill }\mbox{$\spadesuit$}}

\def\orig#1{\bar{#1}}
\newcommand{\ceiling}[1]{\left\lceil #1 \right\rceil}

\def\aa{a}
\def\aao{\bar{a}}
\def\bb{b}
\def\bbo{\bar{b}}
\def\cc{c}
\def\cco{\bar{c}}

\def\xx{x}
\def\yy{b}
\def\zz{c}
\def\tt{t}

\def\origin{{\mbox{\boldmath $0$}}}

\def\form#1#2{\left<#1 | #2 \right>}

\def\prob#1#2{\mbox{\bf Pr}_{#1}\left[ #2 \right]}
\def\expec#1#2{\mbox{\bf E}_{#1}\left[ #2 \right]}

\def\setof#1{\left\{#1  \right\}}
\def\abs#1{\left| #1 \right|}
\def\Span#1{\mbox{{\bf Span}}\left(#1  \right)}
\def\dist#1#2{\mbox{{\bf dist}}\left(#1, #2 \right)}
\def\height#1{\mbox{{\bf height}}\left(#1  \right)}
\def\vs#1#2#3{#1_{#2},\ldots , #1_{#3}}
\def\Reals#1{\StevesR^{#1}}
\def\norm#1{\left\| #1 \right\|}

\begin{document}

\maketitle

\thispagestyle{first} \setcounter{page}{597}

\begin{abstract}\vskip 3mm
Spielman and Teng~[STOC '01] introduced the smoothed analysis of
  algorithms to provide a framework in which one could
  explain the
  success in practice of algorithms and heuristics that could not be
  understood through the traditional worst-case and average-case analyses.
In this talk, we survey some of the smoothed analyses that have
been performed.

\vskip 4.5mm



\end{abstract}

\vskip 12mm

\section{Introduction}\setzero
\vskip-5mm \hspace{5mm }

The most common theoretical approach to understanding
  the behavior of algorithms is worst-case analysis.
In worst-case analysis, one proves a bound on the worst possible
  performance an algorithm can have.
A triumph of the Algorithms community has been the proof that
  many algorithms have good worst-case performance---a strong guarantee
  that is desirable in many applications.
However, there are many algorithms that work exceedingly well in
practice,
  but which are known to perform poorly in the worst-case or
  lack good worst-case analyses.
In an attempt to rectify this discrepancy between theoretical
analysis and
  observed performance, researchers introduced the average-case analysis
  of algorithms.
In average-case analysis, one
  bounds the expected performance of an algorithm on random inputs.
While a proof of good average-case performance provides evidence
that
  an algorithm may perform well in practice, it can rarely be
  understood to explain the good behavior of an algorithm in practice.
A bound on the performance of an algorithm under one distribution
   says little
   about its performance under another distribution, and may
   say little about the inputs that occur in practice.
Smoothed analysis is a hybrid of worst-case and average-case
analyses
  that inherits advantages of both.

In the formulation of smoothed analysis used
in~\cite{SpielmanTengSimplex},
  we measure the maximum over inputs of the expected running time of a
  simplex algorithm under slight random perturbations of
  those inputs.
To see how this measure compares with worst-case and average-case
analysis,
  let $X_{n}$ denote the space of linear-programming problems of
  length $n$ and let
  $\mathcal{T} (x)$ denote the running time of the simplex algorithm
  on input $x$.
Then, the worst-case complexity of the simplex algorithm is the
function
\[
  \mathcal{C}_{worst} (n) = \max_{x \in X_{n}} \mathcal{T} (x),
\]
and the average-case complexity of the algorithm is
\[
  \mathcal{C}_{ave} (n) = \mathbf{E}_{r \in X_{n}} \mathcal{T} (r),
\]
under some suitable distribution on $X_{n}$. In contrast, the
smoothed complexity of the simplex algorithm is the function
\[
  \mathcal{C}_{smooth} (n, \sigma ) = \max_{x} \mathbf{E}_{r \in
X_{n}}
   \mathcal{T} (x + \sigma \norm{x} r),
\]
where $r$ is chosen according to some distribution, such as a
Gaussian. In this case $\sigma \norm{x} r$ is a Gaussian random
vector of
  standard deviation $\sigma \norm{x}$.
We multiply by $\norm{x}$ so that we can relate the magnitude of
the
  perturbation to the magnitude of that which it perturbs.

In the smoothed analysis of algorithms, we measure the expected
performance
  of algorithms under slight random perturbations of worst-case inputs.
More formally, we consider the maximum over inputs of the expected
performance
  of algorithms under slight random perturbations of those inputs.
We then express this expectation as a function of the input size
and
  the magnitude of the perturbation.
While an algorithm with a good worst-case analysis will perform
  well on all inputs, an algorithm with a good smoothed analysis
  will perform well on almost all inputs in every small neighborhood
  of inputs.
Smoothed analysis makes sense for algorithms whose inputs are
  subject to slight amounts of noise in their low-order
  digits, which is typically the case if they are
  derived from measurements of real-world phenomena.
If an algorithm takes such inputs and has a good smoothed
analysis,
  then it is unlikely that it will encounter an
  input on which it performs poorly.
The name ``smoothed analysis'' comes from the observation that if
  one considers the running time of an algorithm as a function from
  inputs to time, then the smoothed complexity of the algorithm is the
  highest peak in the plot of this function after it is convolved
  with a small Gaussian.

In our paper introducing smoothed analysis, we proved that
  the simplex method has polynomial smoothed complexity
  ~\cite{SpielmanTengSimplex}.
The simplex method, which has been the most popular
  method of solving linear programs since the late 1940's,
  is the canonical example of a practically useful algorithm
  that could not be understood theoretically.
While it was known to
  work very well in practice, contrived examples on which
  it performed poorly proved that it had horrible worst-case
  complexity\cite{KleeMinty,Murty,GoldfarbSit,Goldfarb,AvisChvatal,Jeroslow,AmentaZiegler}.
The average-case complexity of the simplex method was proved to
  be polynomial~\cite{Borg82,Borg77,SmaleRand,Haimovich,AdlerKarpShamir,AdlerMegiddo,ToddRand},
  but this result was not considered to explain
  the performance of the algorithm in practice.

\section{The simplex method for linear programming}\setzero
\vskip-5mm \hspace{5mm }

We recall that a linear programming problem can be written in the
form
\begin{eqnarray}
 &  \mbox{maximize} &  \xx^{T} \zz \nonumber  \\
 & \mbox{subject to} & \xx^{T}\aa _{i} \leq b _{i},
  \mbox{ for $1 \leq i \leq n$,}
\end{eqnarray}
where $\zz \in \Reals{d}$, $\aa_{i} \in \Reals{d}$
  and $b_{i} \in \Reals{}$, for $1 \leq i \leq n$
In~\cite{SpielmanTengSimplex}, we bound the smoothed complexity of
a particular
  two-phase simplex method that uses the shadow-vertex pivot rule to
  solve linear programs in this form.

We recall that the constraints of the linear program,
  that $\xx^{T}\aa _{i} \leq b _{i}$, confine $\xx$
  to a (possibly open) polytope, and that the solution
  to the linear program is a vertex of this polytope.
Simplex methods work by first finding some vertex of the polytope,
  and then walking along the 1-faces of the polytope from
  vertex to vertex, improving the objective function at each step.
The pivot rule of a simplex algorithm dictates which vertex
  the algorithm should walk to when it has many to choose from.
The shadow-vertex method is inspired by the simplicity of the
  simplex method in two-dimensions: in two-dimensions, the polytope
  is a polygon and the choice of next vertex is always unique.
To lift this simplicity to higher dimensions, the shadow-vertex
  simplex method considers the orthogonal projection of the polytope
  defined by the constraints onto a two-dimensional space.
The method then walks along the vertices of the polytope
  that are the pre-images of the vertices of the shadow polygon.
By taking the appropriate shadow, it is possible to guarantee that
  the vertex optimizing the objective function will be encountered
  during this walk.
Thus, the running time of the algorithm may be bounded by the
  number of vertices lying on the shadow polygon.
Our first step in proving a bound on this number is a smoothed
  analysis of the number of vertices in a shadow.
For example, we prove the bound:

{\bf Theorem 2.1} (Shadow Size) \it Let $d \geq 3$ and $n > d$.
Let $\zz$ and $\tt$ be independent vectors in $\Reals{d}$,
  and let
  $\vs{\aa }{1}{n}$ be Gaussian
  random vectors
  in $\Reals{d}$ of variance
  $\sigma ^{2} \leq \frac{1}{9 d \log n}$
  centered at points each of norm at most $1$.
Then, the expected number of vertices of the shadow polygon
  formed by the projection of
  $\setof{\xx : \xx^{T} \aa_{i} \leq 1}$
  onto $\Span{\tt ,\zz}$ is at most
\begin{equation}\label{eqn:shadow}
\frac
    {58,888,678 \ n d^{3}
}{
    \sigma ^{6}
}.
\end{equation}
\rm

This bound does not immediately lead to a bound on the running
time of
  a shadow-vertex method as it assumes that $\tt$ and $\zz$ are fixed
  before the $\aa_{i}$s are chosen, while in a simplex method the plane
  on which the shadow is followed depends upon the $\aa_{i}$s.
However, we are able to use the shadow size bound as a black-box
to
  prove for a particular randomized two-phase shadow vertex simplex
  method:

{\bf Theorem 2.2} (Simplex Method) \it Let $d \geq 3$ and $n \geq
d + 1$. Let $\zz \in \Reals{d}$ and $\yy \in \setof{-1,1}^{n}$.
Let
  $\vs{\aa }{1}{n}$ be Gaussian
  random vectors
  in $\Reals{d}$ of variance
  $\sigma ^{2} \leq \frac{1}{9 d \log n}$
  centered at points each of norm at most $1$.
Then the expected number of simplex steps
  taken by the
{two-phase shadow-vertex simplex algorithm}
  to solve the program specified by $\yy$, $\zz$,
  and $\vs{\aa}{1}{n}$ is at most
\[
  (n d / \sigma)^{O (1)},
 \]
where the expectation is over the choice of $\vs{\aa}{1}{n}$
  and the random choices made by the algorithm.
\rm

While the proofs of Theorems~2.1 and~2.2
  are quite involved, we can provide the reader with this intuition
  for Theorem~2.1:
  \textit{after perturbation, most of the vertices of the polytope
  defined by the linear program have an angle bounded away from flat}.
This statement is not completely precise because ``most'' should
be
  interpreted under a measure related to the chance a vertex appears
  in the shadow, as opposed to counting the number of vertices.
Also, there are many ways of measuring high-dimensional
  angles, and different approaches are used in different
  parts of the proof.
However, this intuitive statement tells us that most vertices on
  the shadow polygon should have angle bounded away from flat,
  which means that there cannot be too many of them.

One way in which angles of vertices
  are measured is by the condition
  number of their defining equations.
A vertex of the polytope is given by a set of equations of the
form
\[
  C x = b.
\]
The condition number of $C$ is defined to be
\[
  \kappa (C) =  \norm{C} \norm{C^{-1}},
\]
where we recall that
\[
  \norm{C} = \max_{x} \frac{\norm{C x}}{\norm{x}}
\]
and that
\[
  \norm{C^{-1}} = \min_{x} \frac{\norm{C x}}{\norm{x}}.
\]
The condition number is a measure of the sensitivity of $x$
  to changes in $C$ and $b$,
  and is also a normalized measure of the distance of $C$
  to the set of singular matrices.
For more information on the condition number of a matrix,
  we refer the reader to one of~\cite{GolubVanLoan,TrefethenBau,DemmelBook}.
Condition numbers play a fundamental role in Numerical Analysis,
  which we will now discuss.

\section{Smoothed complexity framework for numerical analysis}\setzero
\vskip-5mm \hspace{5mm }

The condition number of a problem instance is generally defined to
be
  the sensitivity of the output to slight perturbations
  of the problem instance.
In Numerical Analysis, one often bounds the running time of
  an iterative algorithm in terms of the condition number
  of its input.
Classical examples of algorithms subject to such analyses
  include Newton's method for root
  finding and the conjugate gradient method of solving
  systems of linear equations.
For example, the number of iterations taken by the
  method of conjugate gradients is proportional to
  the square root of the condition number.
Similarly, the running times of interior-point methods
  have been bounded
  in terms of condition numbers~\cite{RenegarCond}.

Blum~\cite{BlumSantaFe} suggested that
  a complexity theory of numerical algorithms should be parameterized
  by the condition number of an input in addition to the input size.
Smale~\cite{SmaleComplexity} proposed
  a complexity theory of numerical algorithms in which one:
\begin{enumerate}
\item [1.] \textit{proves a bound on the running time of an algorithm solving
  a problem in terms of its condition number, and then}
\item [2.] \textit{proves that it is unlikely that a random problem
  instance has large condition number.}
\end{enumerate}
This program is analogous to the average-case complexity
  of Theoretical Computer Science and hence shares the same shortcoming
  in modeling practical performance of numerical algorithms.

To better model the inputs that occur in practice,
  we propose replacing step 2 of Smale's program with
\begin{enumerate}
\item [2'.] \textit{prove that for every input instance it is unlikely
  that a slight random perturbation of that instance has
  large condition number.}
\end{enumerate}
That is, we propose to bound the smoothed value of the condition
number. In contrast with the average-case analysis of condition
numbers,
  our analysis can be interpreted as demonstrating
  that if there is a little bit of imprecision or noise in the input,
  then it is unlikely it is ill-conditioned.
The combination of step 2' with step 1 of Smale's program
  provides a simple framework for performing smoothed analysis of
  numerical algorithms whose running time can be bounded by
  the condition number of the input.

\section{Condition numbers of matrices}\setzero
\vskip-5mm \hspace{5mm }

One of the most fundamental condition numbers is the condition
number
  of matrices defined at the end of Section~2.
In his paper, ``The probability that a numerical analysis problem
is difficult'',
  Demmel~\cite{DemmelProb} proved that it is unlikely that
  a Gaussian random matrix centered at the origin
  has large condition number.
Demmel's bounds on the condition number were
   improved by Edelman~\cite{Edelman}.
As bounds on the norm of a random matrix are
  standard, we focus on the norm of the inverse, for which
  Edelman proved:

{\bf Theorem 4.1} (Edelman) \it Let $G$ be a $d$-by-$d$
  matrix of independent Gaussian random variables
  of variance $1$ and mean $0$.
Then,
\[
  \prob{}{\norm{G^{-1}} > t} \leq \frac{\sqrt{d}}{t}.
\]
\rm

We obtain a smoothed analogue of this bound in work
  with Sankar~\cite{SankarSpielmanTeng}.
That is, we show that for every matrix it is unlikely that
  the slight perturbation of that matrix has large condition number.
The key technical statement is:

{\bf Theorem 4.2} (Sankar-Spielman-Teng) \it Let~$\orig{A}$ be an
arbitrary $d$-by-$d$ Real matrix
  and~$A$ a matrix of independent Gaussian random variables centered
  at~$\orig{A}$, each of variance~$\sigma^2$.
Then
  \[ \prob{}{\norm{A^{-1}} > x} < 1.823 \frac{\sqrt d}{x\sigma} \]
\rm

In contrast with the techniques used by Demmel and Edelman,
  the techniques used in the proof of
  Theorem~4.2 are geometric and completely
  elementary.
We now give the reader a taste of these techniques by proving the
  simpler:

{\bf Theorem 4.3 } \it Let~$\orig{A}$ be an arbitrary $d$-by-$d$
Real matrix
  and~$A$ a matrix of independent Gaussian random variables centered
  at~$\orig{A}$, each of variance~$\sigma^2$.
Then,
\[
  \prob{}{\norm{A^{-1}} \geq x }
  \leq
  d^{3/2} / x \sigma.
\]
\rm

The first step of the proof is to relate
  $\norm{A^{-1}}$ to a geometric quantity
  of the vectors in the matrix $A$.
The second step is to bound the probability of a configuration
  under which this geometric quantity is small.
The geometric quantity is given by:

{\bf Definition} \it For $d$ vectors in $\Reals{d}$,  $\vs{a
}{1}{d}$,
  define
\[
  \height{\vs{a }{1}{d}} = \min_{i} \dist{a _{i}}
  {\Span{a _{1},\ldots , \hat{a _{i}}, \ldots ,a _{d}}}.
\]
\rm

{\bf Lemma 4.5} \it For $d$ vectors in $\Reals{d}$,
$\vs{a}{1}{d}$,
\[
  \norm{(\vs{a}{1}{d})^{-1}} \leq \sqrt{d}/ \height{A}.
\]
\rm

{\bf Proof.} Let $\tt $ be a unit vector such that
\[
  \norm{\sum_{i=1}^{d} t_{i} a _{i}} = 1/\norm{(\vs{a}{1}{d})^{-1}}.
\]
Without loss of generality, let $t_{1}$
  be the largest entry of $\tt$
  in absolute value, so $\abs{t_{1}} \geq 1/\sqrt{d}$.
Then, we have
\begin{eqnarray*}
  \norm{a _{1} + \sum _{i=2}^{d} (t_{i}/t_{1}) a _{i}}
   &  \leq &  \sqrt{n}/ \norm{(\vs{a}{1}{d})^{-1}} \\
 \implies
  \dist{a _{1}}{ \Span{\vs{a }{2}{d}}} & \leq &
   \sqrt{d}/ \norm{(\vs{a}{1}{d})^{-1}} .
\ \ \ \ \ \ \qed
\end{eqnarray*}

{\bf Proof of Theorem~4.3.} Let $\vs{a}{1}{d}$ denote the columns
of $A$. Lemma~4.5 tells us that if
  $\norm{A^{-1}} \geq x$, then
  $\height{\vs{a}{1}{d}}$ is less than $\sqrt{d}/x$.
For each $i$,
  the probability that the height of $a _{i}$ above
  $ {\Span{a _{1},\ldots , \hat{a _{i}}, \ldots ,a _{d}}}$
  is less than $\sqrt{d}/ x$
  is at most
\[
    \sqrt{d} / x \sigma.
\]
Thus,
\begin{eqnarray*}
\lefteqn{  \prob{}{\height{\vs{a}{1}{d}} <  \sqrt{d} / x}}\\
& \leq &
  \prob{}{\exists i : \dist{a _{i}}
                      {\Span{a _{1},\ldots , \hat{a _{i}}, \ldots ,a _{n}}}
              <  \sqrt{d}/ x }\\
& < &
   n^{3/2} / x \sigma;
\end{eqnarray*}
so,
\[
  \prob{}{\norm{(\vs{a}{1}{d})^{-1}} > x }
 <
   d^{3/2} / x \sigma. \ \ \ \ \ \  \qed
\]

{\bf Conjecture 1} \it Let~$\orig{A}$ be an arbitrary $d$-by-$d$
Real matrix
  and~$A$ a matrix of independent Gaussian random variables centered
  at~$\orig{A}$, each of variance~$\sigma^2$.
Then
  \[ \prob{}{\norm{A^{-1}} > x} < \frac{\sqrt d}{x\sigma}. \]
\rm

\section{Smoothed condition numbers of linear programs}\setzero
\vskip-5mm \hspace{5mm }

The Perceptron algorithm solves linear programs of the following
  simple form:
\begin{quotation} Given a set of points $\vs{\aa}{1}{n}$, find a
vector $x $ such that $\form{\aa_{i}}{x } > 0$ for all $i$,
 if one exists.
\end{quotation}
One can define the condition number of the Perceptron problem to
be
  the reciprocal of
  the ``wiggle room'' of the input. That is, let
  $S = \{x : \form{\aa_i}{x } > 0, \forall i\}$ and
\[
\nu (\vs{\aa}{1}{n}) = \max_{x \in S}
\left(\min_{i}\frac{\abs{\form{\aa_i}{x}}}{\norm{\aa_i}\norm{x}}\right).
\]
Then, the condition number of Perceptron problem is defined to be
  $1/\nu (\vs{\aa}{1}{n})$.

The Perceptron algorithm works as follows: (1) Initialize $x
=\origin $; (2) Select any $\aa_i$
 such that $\form{\aa_i}{x}\leq 0$ and set $x = x +
\aa_i/\norm{\aa_i}$; (3) while $x \not\in S$, go back to step (2).

Using the following two lemmas, Blum and
Dunagan~\cite{BlumDunagan}
  obtained a smoothed
  analysis of the Perceptron algorithm.

{\bf Theorem 5.1} (Block-Novikoff) \it On input $\vs{\aa}{1}{n}$,
the perceptron algorithm terminates in at most $1/ (\nu
(\vs{\aa}{1}{n}))^{2}$ iterations. \rm

{\bf Theorem 5.2} (Blum-Dunagan) \it Let $\vs{\aa}{1}{n}$ be
Gaussian random vectors in $\Reals{d}$ of variance $\sigma^{2}< 1/
(2d)$ centered at points each of norm at most $1$. Then,
\[
\prob{}{\frac{1}{\nu (\vs{\aa}{1}{n})} > t} \leq \frac{n
d^{1.5}}{\sigma t}\log
  \frac{\sigma t}{d^{1.5}}.
\]
\rm

Setting $t = \frac{n d^{1.5}\log (n/\delta)}{\delta \sigma}$, Blum
and
  Dunagan concluded

{\bf Theorem 5.3} (Blum-Dunagan) \it Let $\vs{\aa}{1}{n}$ be
Gaussian random vectors in $\Reals{d}$ of variance $\sigma^{2}< 1/
(2d)$ centered at points each of norm at most $1$. Then, there
exists a constant $c$ such that
  the probability that the perceptron takes
  more than
$ \frac{c d^{3}n^{2}\log^{2 }
  (n/\delta)}{\delta^{2}\sigma^{2}}$
iterations is at most $\delta$. \rm

In his seminal work,
  Renegar~\cite{RenegarPert,RenegarCond,RenegarFunc} defines the condition
  of a linear program to be the normalized reciprocal of its distance
  to the set of ill-posed linear programs, where
 an ill-posed program is one that can be made both feasible and
  infeasible or bounded and unbounded
   by arbitrarily small changes to its constraints.

Renegar proved the following theorem.

{\bf Theorem 5.4} (Renegar) \it There is an interior point method
such that,
  on input a linear program specified by $(A, \yy , \zz)$
  and an $\epsilon > 0$,
  it will terminate in
\[
O (\sqrt{n+d}\log (\kappa (A,\bb,\cc)/\epsilon )
\]
iterations
  and return either a
  solution within $\epsilon$ of the optimal
  or a certificate that the linear program is infeasible or
  unbounded.
\rm

With Dunagan, we recently proved the following smoothed bound on
  the condition number of a linear program~\cite{DunaganSpielmanTeng}:

{\bf Theorem 5.5} (Dunagan-Spielman-Teng) \it For any $\sigma^{2}<
1/ (nd)$,
  let $A = (\vs{\aa}{1}{n})$ be a set of Gaussian random vectors in
  $\Reals{d}$ of variance $\sigma^{2}$ centered at points
  $\vs{\aao}{1}{n}$, let
  $\bb$ be a Gaussian random vector in $\Reals{d}$ of variance
  $\sigma^{2}$ centered at $\bbo$ and let $\cc$ be a Gaussian random vector
  in $\Reals{n}$ of variance $\sigma^{2}$  centered at $\cco$ such that
  $\sum_{i=1}^{n}\norm{\aao_i}^{2}+ \norm{\bbo}^{2} + \norm{\cco}^{2} \leq
  1$.
Then
\[
\prob{A,\bb,\cc}{C (A,\bb,\cc ) > t} < \frac{2^{14}n^{2}d^{3/2}}
  { \sigma^{2} t}\log^{2}\frac{2^{10}n^{2}d^{3/2}t}{\sigma^{2}},
\]
and hence
\[
\expec{A,\bb,\cc}{\log C (A,\bb,\cc)} \leq  21 + 3\log (nd/\sigma
).
\]
\rm

Combining these two theorem, we have

{\bf Theorem 5.6 }(Smoothed Complexity of Interior Point Methods)
\it Let $\sigma$ and  $(A,\bb,\cc)$
  be as given in Theorem~5.5,
Then, Renegar's interior point method
  solves the linear program
  specified by $(A,\bb,\cc)$ to within precision $\epsilon $
  in expected
\[
O \left(\sqrt{n+d} (21 + 3\log (nd/\sigma\epsilon  ))\right)
\]
iterations. \rm

\section{Two open problems}\setzero
\vskip-5mm \hspace{5mm }

As the norm of the inverse of matrix is such a fundamental
quantity,
  it is natural to ask how the
  norms of the inverses of the $\binom{n}{d}$
  $d$-by-$d$ square sub-matrices
  of a $d$-by-$n$ matrix behave.
Moreover, a crude bound on the probability that many of these are
  large is a dominant term in the analysis of complexity of the simplex
  method in~\cite{SpielmanTengSimplex}.
The bound obtained in that paper is:

{\bf Lemma 6.1} \it Let $\vs{\aa}{1}{n}$ be Gaussian random
vectors in $\Reals{d}$
  of variance $\sigma^{2} \leq 1/9 d \log n$
  centered at points of norm at most $1$.
For $I \in \binom{[n]}{d}$ a $d$-set, let
 $X_{I}$ denote the indicator random variable that is 1
 if
\[
  \norm{\left[\aa_{i} : i \in I\right]^{-1}} \geq
     \frac{\sigma^{2}}
     {8 d^{3/2} n^{7}}.
\]
Then,
\[
 \prob{\vs{\aa}{1}{n}}
      {
  \frac{d}{2} \sum_{I} X_{I}
   <  \ceiling{\frac{n-d-1}{2}} \binom{n}{d-1}
   }
  \geq 1 -  n^{-d} -  n^{-n+d-1} - n^{-2.9d + 1}.
\]
\rm

Clearly, one should be able to prove a much stronger bound than
  that stated here, and thereby improve the bounds on the smoothed
  complexity of the simplex method.

While much is known about the condition numbers of
  random matrices drawn from continuous distributions, much
  less is known about matrices drawn from discrete distributions.
We conjecture:

{\bf Conjecture 2} \it Let $A$ be a $d$-by-$d$ matrix of
independently and uniformly chosen
  $\pm 1$ entries.
Then,
\[
  \prob{}{\norm{A^{-1}} > t}
 \leq \frac{\sqrt{d}}{t} + \alpha ^{n},
\]
for some absolute constant $\alpha  < 1$. \rm

We remark that the case $t = \infty$, when the matrix $A$ is
singular,
  follows from a theorem of Kahn, Komlos and Szemeredi~\cite{KahnKomlosSzemeredi}.

\label{lastpage}

\end{document}